\documentclass[10pt,letterpaper,dvips]{article}

 \usepackage{amssymb,latexsym,amsmath,amsthm}

\usepackage{fullpage}

\usepackage[numbers,sort&compress]{natbib}

\usepackage{url}

 \renewcommand{\div}{\mathop{\mathrm{div}}\nolimits}

\newtheorem{thm*}{Theorem}

\newtheorem{thm}{Theorem}

\newtheorem{dfn}{Definition}

\newtheorem{conj}{Conjecture}

\newtheorem{lemma}{Lemma}

\newtheorem{rem}{Remark}

\newtheorem{ques}{Question}

\newtheorem{cor}{Corollary}

\newtheorem{prop}{Proposition}

\newtheorem{note}{Notation}

\begin{document}
\date{}

\title{ Higher dimensional solutions for a nonuniformly elliptic equation}
\maketitle

\begin{center}
\author{Mostafa Fazly\footnote{The author  is pleased to acknowledge the support of a University of Alberta start-up grant RES0019810.}
}
\\
{\it\small Department of Mathematical and Statistical Sciences, 632 CAB, University of Alberta}\\
{\it\small Edmonton, Alberta, Canada T6G 2G1}\\
{\it\small e-mail: fazly@ualberta.ca}\vspace{1mm}
\end{center}

\vspace{3mm}

\begin{abstract} We prove $m$-dimensional symmetry results, that we call  $m$-Liouville theorems, for  stable and monotone solutions of the following nonuniformly elliptic equation 
\begin{eqnarray*}\label{mainequ}
- \div(\gamma(\mathbf x') \nabla u(\mathbf x)) =\lambda (\mathbf x' ) f(u(\mathbf x)) \ \ \text{for}\ \ \mathbf x=(\mathbf x',\mathbf x'')\in\mathbf{R}^d\times\mathbf{R}^{s}=\mathbf{R}^n,
  \end{eqnarray*}
where $0\le m<n$ and $0<\lambda,\gamma$ are smooth functions and $f\in C^1(\mathbf R)$. The interesting fact is that the decay assumptions on the weight function $\gamma(\mathbf x') $ play the fundamental role in deriving  $m$-Liouville theorems.  We show that under certain assumptions on the sign of the nonlinearity $f$, the above equation satisfies a 0-Liouville theorem. More importantly, we prove that for the double-well potential nonlinearities, i.e. $f(u)=u-u^3$, the above equation satisfies a $(d+1)$-Liouville theorem. This can be considered as a higher dimensional counterpart of the celebrated conjecture of De Giorgi for the Allen-Cahn equation.  The remarkable phenomenon is that the $\tanh$ function that is the profile of monotone and bounded solutions of the Allen-Cahn equation appears towards constructing higher dimensional Liouville theorems. 

\end{abstract}


\section{Introduction}

We study $m$-dimensional symmetry of solutions for the following semilinear elliptic equation with an advection term
\begin{equation} \label{vector}
- \Delta u  + \mathbf a(\mathbf{x}) \cdot \nabla u =  b (\mathbf{x}) f(u) \qquad \mathbf{x}\in  \mathbf{R}^n
\end{equation}
where $\mathbf a: \mathbf{R}^{n} \rightarrow \mathbf{R}^{n}$ is a  smooth vector field, $b \in C ^\infty(\mathbf{R}^n)$ and $f\in C^1(\mathbf{R})$.  Note that if $\mathbf a(\mathbf{x})$ is of gradient form, that is there exists a smooth $c(\mathbf{x}) $ such that $\mathbf a(\mathbf{x}) = \nabla c (\mathbf{x})$, then one can rewrite \eqref{vector} as
\begin{equation} \label{vector-2}
- \Delta u  + \nabla c(\mathbf{x}) \cdot \nabla u = b (\mathbf{x}) f(u)  \qquad \mathbf{x}\in  \mathbf{R}^n.
\end{equation}
If we set $\gamma (\mathbf{x}) = e^{- c (\mathbf{x})}$ and $\lambda (\mathbf{x}) = e^{-c (\mathbf{x})} b(\mathbf{x})$ then we can rewrite \eqref{vector-2} as the following equation in divergence form
\begin{equation}\label{main}
-\div (\gamma (\mathbf x) \nabla u)=\lambda( \mathbf x) f(u)\qquad \mathbf{x}\in  \mathbf{R}^n.
\end{equation}
Therefore, we assume that $ \gamma(\mathbf x)$ and $ \lambda(\mathbf x)$, which we call \emph{weights},  are smooth positive functions (we allow $ \lambda$ to be zero at say a point) and which satisfy various growth conditions at infinity.  Note that the assumption $\gamma(\mathbf x)>0$ implies that the operator $\div(\gamma(\mathbf x) \nabla \cdot)$ is a nonuniformly elliptic operator.  

\begin{note} Throughout the paper we use the following notations. 
 \begin{itemize}
 \item  The weight functions $\lambda$ and $\gamma$ are only functions of $d$-variables meaning that $\gamma(\mathbf x)=\gamma(\mathbf x')$ and $\lambda(\mathbf x)=\lambda(\mathbf x')$  where $\mathbf x=(\mathbf x',\mathbf x'')\in \mathbf{R}^d\times\mathbf{R}^{s}=\mathbf{R}^n $ for $n=d+s$. Another representation for $x$ in $n$ dimensional space is $\mathbf x=(\mathbf x''',x_n)\in \mathbf{R}^{n-1}\times\mathbf{R} $. 

\item The following class of nonlinearities appears in our results,
\begin{equation*}
\mathcal{G}:=\left\{  g:\mathbf R^+\to \mathbf R^+,  \ \text{is nondecreasing and}  \ \int_{1}^{\infty} \frac{1}{rg(r)} dr =\infty\right\}.  
\end{equation*}
Note that $\mathcal{G}$ is not empty, e.g., $g(r)=\log (1+r)$ is in $\mathcal G$. 
 \end {itemize}
  \end{note}
  The class $\mathcal{G}$ of nonlinearities was  defined by Karp in \cite{k1,k2} and was used by Moschini in \cite{mos}.

\begin{dfn}
We say that (\ref{main}) satisfies $m$-Liouville theorem if for certain $\lambda$ and $\gamma$ solutions of (\ref{main}) are $m$-dimensional for $0\le m<n$, i.e., they exactly depend  on $m$ variables. Similarly, we say that (\ref{main}) satisfies at most $m$-Liouville theorem if solutions of (\ref{main}) are at most $m$-dimensional for $0\le m<n$, i.e., they depend  on at most $m$ variables.
 \end{dfn}
\begin{dfn} We call a classical solution $u$ of (\ref{main}) to be
\begin{enumerate}
\item [(i)] asymptotically convergent if 
\begin{equation}\label{asymp}
\lim_{x_n\to\pm\infty } u(\mathbf{x}''',x_n)\to \pm 1 \ \ \ \text{for all}\ \ \mathbf{x}'''\in \mathbf{R}^{n-1}.
\end{equation}
 If this limit is uniform then we call it uniformly asymptotically convergent.
 \item [(ii)]  monotone if $\partial_{x_n} u(\mathbf x)>0$ for all $\mathbf x\in\mathbf{R}^n$.
\item [(iii)] pointwise stable if there exists a function $0<v$ that satisfies the linearized equation $$-\div(\gamma(\mathbf x') \nabla v) = \lambda(\mathbf x') f'(u) v \ \ \ \text{for all}\ \ \mathbf{x}\in \mathbf{R}^{n}$$
\item[(iv)] stable if for all $\psi\in C_c^1(\mathbf{R}^n)$ the following inequality holds,
  \begin{equation}\label{stability}
 \int_{\mathbf{R}^n} \lambda(\mathbf x') f'(u) \psi^2 d \mathbf{x} \le \int_{\mathbf{R}^n} \gamma (\mathbf x') |\nabla \psi|^2 d \mathbf{x}.
 \end{equation}
\end{enumerate}
 \end{dfn}
 Note that by taking derivative of  (\ref{main}) with respect to $x_n$ monotonicity implies pointwise stability and multiplying  $\frac{\psi^2}{v}$ and doing integration by parts one can see that  pointwise stability implies stability as it is given in \cite{cf}.    The equation (\ref{vector}) is a perturbation of the following semilinear elliptic equation
\begin{equation} \label{semi}
- \Delta u   =  f(u) \qquad \text{in} \ \  \mathbf{R}^n.
\end{equation}   
 When $m=0$ and $1$, $m$-Liouville theorems for (\ref{semi}) known as Liouville theorems and one dimensional symmetry results, respectively, are extensively studied in the literature  \cite{aac,ac,bbg,bcn,bhm,dkw,df,f,gg1,gg2,mm,mos,s}.     The most well-known 1-Liouville theorem is  the following conjecture of De Giorgi in 1978. 
\begin{conj}\label{conj} Suppose that $u:\mathbf{R}^n\to [-1,1]$ is a classical monotone solution of  (\ref{vector}) for $\mathbf a=0$, $b=1$ and $f(u)=u-u^3$.  Then for at least $n\le 8$ equation (\ref{main}) satisfies 1-Liouville theorem.
\end{conj}
From the definition of 1-Liouville theorem, $u$ depends only on one variable and therefore it has to be of the form 
\begin{equation}\label{tanh}
u(\mathbf x)=\tanh \left(\frac{\mathbf x\cdot \mathbf y - c}{\sqrt 2}  \right) \  \ \text{for all} \ \mathbf x\in\mathbf{R}^n
\end{equation}
for some $c\in\mathbf{R}$ and some $\mathbf y \in \mathbf{R}^n$ where $|\mathbf y|=1$ and $y_n>0$. Note that the function $w(t)=\tanh (t/\sqrt 2)$ is the unique solution up to translation of the following ordinary differential equation,
$$-w''=w-w^3, \ w'>0, \ w(\pm\infty)=\pm 1.$$
In 1997,  Ghoussoub and Gui \cite{gg2}  proved the De Giorgi's conjecture for $n= 2$. They used a linear 0-Liouville theorem for the ratio $\sigma := \frac{\partial u}{\partial_{x_1}}/ \frac{\partial u}{\partial_{x_2}}$ developed by Berestycki, Caffarelli and Nirenberg in \cite{bcn} for the study of symmetry properties of positive solutions of semilinear elliptic equations in half spaces. Unfortunately, it is not known whether or not this 0-Liouville theorem is optimal, see Proposition \ref{liouville} and what follows shortly after. 

Ambrosio and Cabr\'{e} \cite{ac} and later in a joint work with Alberti \cite{aac} extended these results up to dimension $n= 3$. The De Giorgi's conjecture for higher dimensions is still open. However,  Ghoussoub and Gui showed in \cite{gg1} that the conjecture is true for $n = 4$ or $n= 5$ for a special class of solutions that satisfy an antisymmetry condition. In 2003, Savin \cite{s} assuming the additional natural hypothesis
 \begin{equation}\label{asym}
\lim_{x_n\to\pm\infty } u(\mathbf{x}''',x_n)\to \pm 1 \ \ \text{for all} \ \ \mathbf{x}'''\in \mathbf{R}^{n-1} ,
\end{equation}
  proved that the conjecture is true in dimension $n\le 8$. The proof is nonvariational and it uses the sliding method for a special family of radially symmetric functions. Finally in 2008,  del Pino-Kowalczyk-Wei in \cite{dkw} gave a counterexample to De Giorgi's conjecture in dimension $n \ge 9$ which has long been believed to exist.  Very recently in \cite{fg}, Ghoussoub and the author gave an extension of the De Giorgi's conjecture to elliptic systems and provided an affirmative answer to the conjecture in lower dimensions.    See also \cite{ali,abg} for more information about the elliptic systems.  
 
Under a much stronger assumption that the limits in \eqref{asym} are uniform in $\mathbf{x}'''$, the conjecture is known as  {\it Gibbons' conjecture}. This conjecture was first proved for $n\le  3$ by Ghoussoub and Gui in \cite{gg2} and then for all dimensions independently  with different methods  by Barlow, Bass and Gui \cite{bbg}, Berestycki, Hamel and Monneau \cite{bhm} and Farina \cite{f}. We also refer interested readers to  \cite{am,aim,sv,fazly} and references therein for some  results regarding the weighted Allen-Cahn equation and system.

 In this article, we attempt to partially answer this question:
\begin{ques}\label{question1}
Under what conditions on $\lambda$ and $\gamma$, $m$-Liouville theorems hold for (\ref{main}) when $0\le m<n$? 
\end{ques}
In other words, we are interested to explore  how a lower order perturbation of  equation (\ref{semi}) would change the behaviour of the solutions. In \cite{cf}, we have proved 0-Liouville theorem in certain dimensions for (\ref{main}) with specific nonlinearities  $f(u) = e^u$, $u^p$ where $p > 1$ and $-u^{-p}$ where $p > 0$ known as the Gelfand, Lane-Emden and negative exponent nonlinearities, respectively. Note that these nonlinearities are non sign changing functions. To prove 0-Liouville theorem, we assumed either $ \nabla \gamma(\mathbf{x}) \cdot \mathbf{x} \le 0$ (i.e. $\nabla c (\mathbf{x}) \cdot \mathbf{x} \ge 0$) or $|\nabla \gamma(\mathbf x) |\le C \lambda(\mathbf x)$ (i.e. $|\nabla c|$ bounded). In this note we prove 0-Liouville theorem for  (\ref{main}) with a general nonlinearity $f\in C^1(\mathbf{R})$ as well as $m$-Liouville theorems for (\ref{main}) when $m\ge 1$ under certain  conditions on $\lambda$ and $\gamma$.     To prove  higher dimenional Liouville theorems, which are more challenging problems, we apply a standard linear 0-Liouville theorem given in \cite{bcn,gg2}.   

The organization of the paper is as follows. In Section \ref{secmain},  we state the main results of the paper and in particular the applications of the main results for the nonuniformly elliptic Allen-Cahn equation.    In Section \ref{seclin},  we provide a linear 0-Liouville theorem and a geometry Poincar\'{e} inequality that are the essential tools in our proofs.  Finally in Section \ref{secm}, we provide $m$-Liouville theorems and the proof of main results.

\section{Main results and related backgrounds}\label{secmain}

  As shown by Gilbarg and Serrin in \cite{gs} (see P. 324) a 0-Liouville theorem holds for bounded solutions of the  linear equation 
\begin{equation}\label{gilbarg}
-\Delta u+\mathbf a(\mathbf x)\cdot \nabla u=0 \ \ \ \text{in $\mathbf{R}^n$}
\end{equation}
where $n\ge 2$ and $\mathbf a(\mathbf x)=O(|\mathbf x|^{-1})$.  The validity or the failure of $0$-Liouville theorems for this equation under appropriate conditions have been also studied in \cite{sssz,stamp}. If we replace the equality with the inequality $\ge$ in (\ref{gilbarg}), then it is strightforward to construct nonconstant bounded solutions satisfying specific $\mathbf a(\mathbf x)=O(|\mathbf x|^{-1})$. This implies a natural question that under what assumptions on $\mathbf a$, $b$ and solutions  one can prove a 0-Liouville theorem for the nonlinear case, (\ref{vector}), with a general nonlinearity $f\ge 0$. In what follows, we prove a 0-Liouville theorem   for bounded stable solutions of  (\ref{vector}). 

\begin{thm}\label{liou3}
Let $u$ be a bounded pointwise stable solution for (\ref{main})  and let either $0\le  f(t)$ or $t f(t)\le 0$ for all $t$ in the range of $u$. If either  
\begin{equation}\label{liou3g}
\int_{B_{2R}}  \gamma(\mathbf x') d\mathbf x' \le k  g(R) \ \ \text{and} \ \  n\le d+4,
\end{equation}
or
 \begin{equation}\label{liou3rg}
\int_{B_{2R}}  \gamma(\mathbf x') d\mathbf x' \le k R  g(R) \ \ \text{and} \ \  n\le d+3,
\end{equation}
where  $g\in \mathcal{G}$ and  $k$ is a constant independent of $R$.  Then, (\ref{main}) satisfies 0-Liouville theorem.
\end{thm}

The proof of the theorem  is strongly motiveated by the methods and ideas developed  by Dupaigne-Farina in \cite{df}, where they examined the advection free equation that is (\ref{vector}) when $\mathbf a=0$ and $b=1$.   Note that the double-well potential nonlinearity $f(t)=t-t^3$ for $t\in [-1,1]$ does not satisfy neither  $0\le  f(t)$ nor $t f(t)\le 0$.  Therefore, in what follows we focus on this type nonlinearity. Berestycki, Hamel and Monneau, Theorem 2 in \cite{bhm},  have shown that a 1-Liouville theorem holds for uniformly asymptotically convergent solutions of (\ref{vector}) under the assumption that $\mathbf a$ is a constant vector, $b(\mathbf x)=b(x_n)$ is bounded and $f$ is Lipschitz continuous on [-1,1]  satisfying 
\begin{enumerate}
\item[(P)] $f(\pm 1)=0$ and  there exists $\delta>0$ such that $f$ is non-increasing on $[-1,-1+\delta]$ and on $[1-\delta,1]$.  
\end{enumerate}
However, a counterexample given by Bonnet-Hamel in \cite{bh} shows that this result no longer holds  if we drop the ''uniformly" assumption. In other words,  they constructed a two dimensional monotone and asymptotically convergent solution such that for $\alpha\in(0,\frac{\pi}{2}]$
\begin{eqnarray*}
u(t \cos\theta, t\sin\theta)&\to& -1 \ \ \text{as} \ \ t\to\infty \ \ \text{for} \ \ -\frac{\pi}{2}-\alpha<\theta< -\frac{\pi}{2}+\alpha\\
u(t \cos\theta, t\sin\theta)&\to& 1 \ \ \text{as} \ \ t\to\infty \ \ \text{for} \ \ -\frac{\pi}{2}+\alpha<\theta< \frac{3\pi}{2}-\alpha
\end{eqnarray*}
when $u$ is a solution of the following equation
  \begin{equation}\label{firstorderC}
  -\Delta u+k \partial_{x_2} u=f(u) \ \ \ \text{in $\mathbf{R}^2$}
  \end{equation}
 where  $k$ is just a constant and for some particular $f$ that satisfies  (P). The level sets of such a solution are parallel lines and cannot be one dimensional. Therefore, De Giorgi's conjecture does not hold for (\ref{firstorderC}). Note that this is a sharp result, since when $k=0$, it follows from the result of Ghoussoub and Gui \cite{gg2} that (\ref{firstorderC}) satisfies a 1-Liouville theorem.

Moreover, Berestycki, Hamel and Monneau, Theorem 3 in \cite{bhm}, have proved   that the 1-Liouville theorem no longer holds for (\ref{vector}) if $a$  is a non constant vector, even for uniformly asymptotically convergent solutions. More precisely, they proved that the following equation in two dimensions
\begin{equation}\label{firstorder}
-\Delta u+a(x_1) \partial_{x_1} u= f(u) \ \ \ \text{in $\mathbf{R}^2$}
\end{equation}
admits both a solution depending on only $x_2$ and infinitely many nonplanar solutions, that is, solutions whose level sets are not parallel.  The construction of nonplanar solutions is very technical and relies on the subsolution-supresolution method. As a conclusion, the Gibbons' conjecture (and therefore De Giorgi's conjecture) cannot be extended to (\ref{firstorder}) that is in dimension two.

In what follows, we provide a higher dimensional Liouville theorem for solutions of (\ref{main}) under certain decay assumptions on $\gamma$ and $\lambda$, and in a particular case this can be applied to prove higher dimensional Liouville theorems for (\ref{firstorder}).  

\begin{thm}\label{liou2} Assume that $f\in C^1([-1,1])$ and $F(t)\le \min\{F(-1),F(1)\}$ for all $t\in (-1,1)$, where $F'=f$. Let  $u$ be a monotone and asymptotically convergent solution of (\ref{main}). Moreover, suppose that there exists a positive constant $k$ such that $|\nabla \gamma (\mathbf x')|\le k \gamma(\mathbf x')$ and $\lambda(\mathbf x') \le k \gamma(\mathbf x')$  for any $\mathbf x'$ outside a compact set in $\mathbf{R}^d$ and  either 
\begin{equation}\label{}
\int_{B_{R}}  \gamma(\mathbf x') d\mathbf x' \le k  g(R) \ \ \text{and} \ \  n\le d+3,
\end{equation}
or
 \begin{equation}\label{}
\int_{B_{R}}  \gamma(\mathbf x') d\mathbf x' \le k R  g(R) \ \ \text{and} \ \  n\le d+2,
\end{equation}
where $g\in \mathcal{G}$ and  $k$ is a constant independent of $R$.  Then, (\ref{main}) satisfies at most $(d+1)$-Liouville theorem.
\end{thm}

\begin{cor}\label{3+epsilon}
Assume that $d=1$, $a\in L^{\infty}(\mathbf R)$ and 
\begin{equation}\label{growth}
\int_{-R}^{R} e^{  -\int_{0}^{x_1} a(t) dt }dx_1 \le k R^{1-\epsilon} g(R)
\end{equation}
where $k=k(n,a,g)$ is a constant independent of $R$,   $0\le \epsilon \le 1$ and any $g\in \mathcal{G}$. Then  monotone and asymptotically convergent solutions of $$-\Delta u+a(x_1) \partial_{x_1} u=u-u^3\ \ \ \text{in} \ \ \mathbf{R} \times\mathbf{R}^{n-1}$$ satisfy at most 2-Liouville theorm for $n\le 3+\epsilon$. 
 \end{cor}

 In particular, this shows that monotone and asymptotically convergent solutions of (\ref{firstorder}) on $\mathbf{R}^n$ and up to dimension $n\le 4$ are at most two dimensional provided $a\in L^{\infty}(\mathbf{R})$ and  
 \begin{equation}\label{consgrowthgamin}
\lim_{R\to\infty} \int_{-R}^{R} \gamma(x_1) dx_1 < \infty \ \  \text{or equivalently} \ \  \lim_{R\to\infty} \int_{-R}^{R}  e^{-\int_{0}^{x_1} a(s)ds}dx_1 < \infty.
\end{equation}
Note that $a(x_1)\equiv k$ where $k$ is just a constant does not satisfy this condition. However, either $a(x_1)=\frac{2x_1}{1+x_1^2}$ or $a(x_1)=t \tanh x_1+s$ for any $t>|s|$ can be chosen to fulfill the assumption (\ref{consgrowthgamin}).  Note also that  the double-well potential $f(t)=t-t^3$  and therefore $F(t)=-\frac{1}{4}(1-t^2)^2$ satisfies the assumptions of Theorem \ref{liou2}. For $\lambda=\gamma=1$ this result is given by Ambrosio-Cabr\'{e}  in \cite{ac} and Ghoussoub-Gui in \cite{gg2}.

The remarkable phenomenon is that according to the De Giorgi's conjecture  monotone  and asymptotically convergent solutions of the Allen-Cahn equation, i.e., (\ref{semi}) with $f(u)=u-u^3$, are one dimensional solutions up to dimension eight and the profile solution is the $\tanh$ function. Now, if we perturb the Allen-Cahn equation by $\tanh$ function that is 
\begin{equation}\label{tanhsemi}
 -\Delta u+\tanh(x_1) \partial_{x_1} u =u-u^3 \ \ \text{in} \ \ \mathbf{R}^n 
\end{equation}
then according to Theorem \ref{liou2} the monotone  and asymptotically convergent  solutions are at most two dimensional up to dimension four. Similarly, higher dimensional Liouville theorems can be constructed as following.

\begin{cor}
Assume that $d=2$, $a_1,a_2\in L^{\infty}(\mathbf R)$ and 
\begin{equation}\label{growthd2}
\int_{x_1^2+x_2^2\le R^2} e^{  -\int_{0}^{x_1} a_1(t) dt }e^{  -\int_{0}^{x_2} a_2(t) dt }dx_1dx_2 \le k R^{2-\epsilon} g(R)
\end{equation}
 where $k=k(n,a_1,a_2,g)$ is a constant independent of $R$,   $0\le \epsilon \le 2$ and any $g\in \mathcal{G}$. Then  monotone and asymptotically convergent solutions of $$-\Delta u+a_1(x_1) \partial_{x_1} u+a_2(x_2) \partial_{x_2} u =u-u^3\ \ \text{in} \ \ \mathbf{R}^2 \times\mathbf{R}^{n-2}$$ satisfy at most 3-Liouville theorem for $n\le 3+\epsilon$. 
 \end{cor}

 Following ideas given in \cite{sz,fsv,sv} we provide  a geometric Poincar\'{e} inequality  for stable solutions of (\ref{main}).  The interesting point is that both the weight function $\lambda$ and the nonlinearity $f$ in (\ref{main})  do not appear in this  geometric Poincar\'{e} inequality. However, the weight function $\gamma$ in (\ref{main}) appears as a weight function for both sides of the inequality.  

\begin{thm}\label{poincare}
Let $u$ be a stable solution of (\ref{main}). Then the following inequality holds for any $\phi\in C_c^1(\mathbf R^n)$,
\begin{eqnarray}\label{poincarein}
&&\int_{\mathbf x'\in\mathbf R^d} \gamma(\mathbf x') \int_{\mathbf x''\in \mathbf R^s \cap \{ \nabla_{\mathbf x''} u \neq0 \}    } \phi^2 \left(  |\nabla_{\mathbf x''} u|^2 \mathcal K ^2 +|\nabla_T |\nabla _{\mathbf x''}u | |^2  \right)  d \mathbf x'' d \mathbf x'  
\\ & &+\nonumber  \int_{\mathbf R^n}  \gamma(\mathbf x')  \phi^2 \mathcal S  \le  \int_{\mathbf R^n}   \gamma(\mathbf x') |\nabla_{\mathbf x''} u|^2 |\nabla \phi|^2
\end{eqnarray}
where $\nabla_T$ denotes the orthogonal projection of the gradient along this level set and
\begin{equation}\label{S}
\mathcal S:=\sum_{j=d+1}^{n}\sum_{i=1}^{d}  |\partial_i\partial_j u|^2  - |  \nabla_{\mathbf x'} |\nabla_{\mathbf x''} u|  |^2
\end{equation}
  and $\mathcal K$ is the full curvature defined by $$\mathcal K( \mathbf x)=\sqrt{ \sum_{j=1}^{s-1} { \mathcal{\kappa}_j  ( \mathbf x) ^2}   }$$ when  $ \mathcal{\kappa}_j$ are the principal curvatures of the level set of $u$ at $\mathbf x$.
\end{thm}

\begin{rem}
The function $\mathcal S$ given in (\ref{S}) is nonnegative.  This can be seen by taking the gradient of $ |  \nabla_{\mathbf x''} u|$ with respect to $\mathbf x'$ and then applying the Cauchy inequality for the points that $|  \nabla_{\mathbf x''} u|\neq0$. 
\end{rem}

In this context and for the case of $\gamma=\lambda=1$, this type of geometric Poincar\'{e} inequality was introduced by Sternberg and Zumbrun in \cite{sz} to study semilinear phase transitions problems. Later on and for the first time, Farina, Sciunzi and Valdinoci in \cite{fsv} used and extended the inquality to prove very interesting results related to the De Giorgi's conjecture.  Then Cabr\'{e} used it (see Proposition 2.2 in \cite{c}) to prove the boundedness of extremal solutions of semilinear elliptic equations with Dirichlet boundary conditions on a convex domain up to dimension four.  Similar inequalities are proved by Savin and Valdinoci  in \cite{sv} for (\ref{main}) when $\gamma=1$.  Recently in \cite{fg}, Ghoussoub and the author extended this inequality to elliptic systems and used it to prove  De Giorgi type results for systems.

\section{Linear 0-Liouville Theorem and a geometric Poincar\'{e} inequality}\label{seclin}

We start this section with the following linear 0-Liouville theorem that is given by Berestycki-Caffarelli-Nirenberg \cite{bcn} and Ghoussoub-Gui  \cite{gg2} for bounded $h\sigma$ and then improved by Ambrosio-Cabr\'{e} \cite{ac} and Moschini \cite{mos}.

\begin{prop} \label{liouville}
Let $0<h\in L^\infty_{loc}(\mathbf{R}^n)$ and $\sigma\in H^1_{loc}(\mathbf{R}^n)$. If $\sigma$ satisfies the following differential inequality  
\begin{equation}\label{harmonic}
\sigma \div ( h(\mathbf x) \nabla \sigma)\ge0 \ \ \ \text{in}\ \ \mathbf{R}^n,
\end{equation}
 such that for any $R>1$,  
\begin{equation}\label{lioubound}
\int_{B_{2R}\setminus B_{R}} h(\mathbf x) \sigma^2 \le C R^2 g(R), 
\end{equation}
where $g\in \mathcal{G} $.  Then $\sigma$ is constant.
\end{prop} 
Note that in two dimensions Proposition \ref{liouville} is sharp in the sense that the following example 
$$h:\equiv1\ \ \text{and for} \ R_0>e^{3/4} \  \text{set} \ \sigma(r):=
 \left\{\begin{array}{lcl}
  \log R_0+\frac{r^2}{R_0^2}-\frac{r^4}{4R_0^4}-\frac{3}{4} \ \ &\text{for}&\ \ \ r<R_0,\\ 
\log r \ \ &\text{for}& \ \ \ r\ge R_0,
\end{array}\right.
$$
given in \cite{mos} (Remark 5.4) shows that this proposition does not hold if $g(R) =\log^2(R)$. Straightforward calculations show that  $\log^2(1+r)$ is not in the class $\mathcal{G}$, however $\log(1+r)$ belongs to $\mathcal{G}$.  

Ambrosio and Cabr\'{e} in \cite{ac} and later on with Alberti in \cite{aac} proved the the following energy estimate holds in any dimension regarding the De Giorgi's conjecture 
\begin{equation}\label{energy}
\int_{B_R} |\nabla u|^2 \le C R^{n-1}.
\end{equation}
Then applying Proposition \ref{liouville} when $g=1$ and equating the right hand sides of (\ref{energy}) and (\ref{lioubound}) they gave a positive answer to Conjecture \ref{conj} in  three dimensions. Now, comparing (\ref{energy}) and (\ref{lioubound})  in any dimensions for the choice of $g(R)=R^{n-3}$,  one sees that the right-hand side of these integral estimates are the same.  Therefore, potentially the function $g(R)=R^{n-3}$ can play an important role in solving  Conjecture \ref{conj} in dimensions $4\le n\le 8$. 

In \cite{mos} and as Remark 5.5,  it's been asked to prove or disprove Proposition \ref{liouville} when $g(R)=R^{n-3}$ and $4\le n\le 8$.  This is a very interesting question because this choice of function $g(R)$ does not belong to the class of $\mathcal{G}$ and also for given function $g(R)$ and in dimensions $n \ge 9$, Ghoussoub and Gui constructed a counterexample for  Proposition \ref{liouville},  see Proposition 2.6 in \cite{gg2}. Their counterexample is very well-constructed and satisfies $\sigma \div ( h(\mathbf x) \nabla \sigma)=0$.   Here we give an elementary example that shows for the subsolution case (inequality $\ge$ holds in (\ref{harmonic})) Proposition \ref{liouville}  does not hold   when $g(R)=R^{n-3}$ and $4\le n\le 8$.

 \begin{rem} Let $n\ge 4$, $h(\mathbf x)=(1+|\mathbf x|^2)^{-\frac{2n-5}{2}}$ and $\sigma(\mathbf x)=(1+|\mathbf x|^2)^\frac{n-3}{2}$. The functions $h$ and $\sigma$  are smooth functions and $0<h\in L^\infty(\mathbf R^n)$. By a simple calculation one can see  that (\ref{harmonic}) holds and moreover $$\int_{B_R}  h(\mathbf x)\sigma^2\le R^{n-1}=R^2 g(R)$$
  where $g(R)=R^{n-3}$. Therefore, $h$ and $\sigma$ satisfy the assumptions of Proposition \ref{liouville}. But $\sigma$ is not a constant even though $h\sigma^2\in L^\infty(\mathbf R^n)$.  This means that to prove Conjecture \ref{conj} in dimensions $4\le n\le 8$ via using (\ref{energy}) and (\ref{lioubound}) when $g(R)=R^{n-3}$,  a counterpart of  Proposition \ref{liouville} is needed that assumes equality in (\ref{harmonic}) and allows a wider class of functions in $\mathcal{G}$.
  \end{rem}
  
  For the rest of this section, we provide a proof for the geometric Poincar\'{e} inequality (\ref{poincarein}).  
\\
\\
\noindent\textbf{Proof of Theorem \ref{poincare}:}   Let $u$ be a stable solution of (\ref{main}). Test the stability inequality (\ref{stability}) with $\psi=|\nabla_{\mathbf x''} u| \phi$ where $\phi\in C_c^1(\mathbf R^n)$ is a test function to get 
  \begin{equation}\label{IJ}
 I:= \int_{\mathbf{R}^n} \lambda(\mathbf x') f'(u) |\nabla_{\mathbf x''} u|^2 \phi^2 \le \int_{\mathbf{R}^n} \gamma (\mathbf x') |\nabla      \left( |\nabla_{\mathbf x''}u |   \phi\right)  | ^2  =:J
 \end{equation}
In what follows we simplify $I$ and $J$. Let's start with $I$.  
\begin{eqnarray}\label{I}
\nonumber I&=& \int_{\mathbf{R}^n} \lambda(\mathbf x') f'(u)  \nabla_{\mathbf x''} u\cdot \nabla_{\mathbf x''} u \ \phi^2= \int_{\mathbf{R}^n} \nabla_{\mathbf x''} \left(  \lambda(\mathbf x')  f(u)   \right)\cdot  \nabla_{\mathbf x''} u \ \phi^2 
\nonumber 
\\&=& - \int_{\mathbf{R}^n} \nabla_{\mathbf x''} \left(   \div (  \gamma(\mathbf x') \nabla u   )   \right)   \cdot  \nabla_{\mathbf x''} u \ \phi^2 = - \int_{\mathbf{R}^n} \nabla_{\mathbf x''}  \left[   \gamma(\mathbf x')  \Delta u + \nabla_{\mathbf x'}\gamma(\mathbf x') \cdot \nabla_{\mathbf x'} u    \right]\cdot  \nabla_{\mathbf x''} u   \ \phi^2 
\nonumber 
\\&=&  - \sum_{i=d+1}^{n}    \int_{\mathbf{R}^n}  \gamma(\mathbf x')  \partial_i u \Delta (\partial_i u) \ \phi^2  -  \sum_{i=d+1}^{n}  \sum_{j=1}^{d} \int_{\mathbf{R}^n} \partial_{j}\gamma(\mathbf x') \partial_{i} \partial_{j} u\partial_iu \ \phi^2 
\nonumber 
\\&=&    \sum_{i=d+1}^{n} \sum_{j=1}^{n}   \int_{\mathbf{R}^n}  \gamma(\mathbf x')  \partial_i u  \partial_j (\partial_i u) \partial_j (\phi^2) +   \sum_{i=d+1}^{n} \sum_{j=1}^{n}   \int_{\mathbf{R}^n}  \gamma(\mathbf x')   |\partial_j \partial_i u|^2 \ \phi^2  
\nonumber \\&=& \frac{1}{2} \sum_{i=d+1}^{n}  \int_{\mathbf{R}^n}  \gamma(\mathbf x')   \nabla (  \partial_i u  )^2 \cdot \nabla \phi^2    +   \int_{\mathbf{R}^n}  \gamma(\mathbf x')  | D^2_{\mathbf x''} u |^2 \ \phi^2 +\sum_{i=d+1}^{n} \sum_{j=1}^{d}   \int_{\mathbf{R}^n}  \gamma(\mathbf x')   |\partial_j \partial_i u|^2 \ \phi^2
\end{eqnarray}
Note that fortunately the term that includes the gradient of $\gamma$ cancels out in the fourth line of calculations where we have used integration by parts. Now we simplify the integral term given as $J$. First note that for $\psi=|\nabla_{\mathbf x''} u| \phi$ we have
\begin{eqnarray*}\label{partJ}
|\nabla \psi|^2 &=&   \left| \nabla |  \nabla_{\mathbf x''} u | \right|^2 \phi^2 + |  \nabla_{\mathbf x''} u|^2 |\nabla \phi|^2 +\frac{1}{2} \nabla \phi^2\cdot \nabla  |  \nabla_{\mathbf x''} u |^2
\\&=&\left| \nabla_{\mathbf x'} |  \nabla_{\mathbf x''} u | \right|^2 \phi^2 + \left| \nabla_{\mathbf x''} |  \nabla_{\mathbf x''} u | \right|^2 \phi^2 + |  \nabla_{\mathbf x''} u|^2 |\nabla \phi|^2 +\frac{1}{2} \nabla \phi^2\cdot \nabla  |  \nabla_{\mathbf x''} u |^2.
\end{eqnarray*}
Therefore, 
\begin{eqnarray}\label{J}
 \nonumber \int_{\mathbf{R}^n} \gamma (\mathbf x') |\nabla      \left( |\nabla_{\mathbf x''}u |   \phi\right)  | ^2 & =&  \int_{\mathbf{R}^n}  \gamma (\mathbf x')   \left| \nabla_{\mathbf x'} |  \nabla_{\mathbf x''} u | \right|^2 \phi^2 + \int_{\mathbf{R}^n}  \gamma (\mathbf x')  \left| \nabla_{\mathbf x''} |  \nabla_{\mathbf x''} u | \right|^2 \phi^2\\&& +\int_{\mathbf{R}^n}  \gamma (\mathbf x')   |  \nabla_{\mathbf x''} u|^2 |\nabla \phi|^2 +\frac{1}{2} \int_{\mathbf{R}^n}  \gamma (\mathbf x')  \nabla \phi^2\cdot \nabla  |  \nabla_{\mathbf x''} u |^2
\end{eqnarray}
The first term in the right-hand side of  (\ref{I}) and the last term in the right-hand side of (\ref{J}) are the same.  Substituting  (\ref{I}) and (\ref{J}) in (\ref{IJ}) we get 
\begin{eqnarray}\label{stend}
\int_{\mathbf{R}^n}  \gamma(\mathbf x') \left(  | D^2_{\mathbf x''} u |^2 - \left| \nabla_{\mathbf x''} |  \nabla_{\mathbf x''} u | \right|^2 \right) \ \phi^2 +  \int_{\mathbf R^n}  \gamma(\mathbf x')  \mathcal S   \phi^2 \le \int_{\mathbf{R}^n}  \gamma (\mathbf x')   |  \nabla_{\mathbf x''} u|^2 |\nabla \phi|^2
\end{eqnarray}
 According to formula (2.1) given in  \cite{sz}, the following geometric identity between the tangential gradients and curvatures holds. For any $w \in C^2(\Omega)$ where $\Omega $ is an open set in $\mathbf R^s$
 \begin{eqnarray}\label{identity} \sum_{i=1}^{s} |\nabla \partial_k w|^2-|\nabla|\nabla w||^2=
\left\{
                      \begin{array}{ll}
                       |\nabla w|^2 (\sum_{i=1}^{s-1} \mathcal{\kappa}_l^2) +|\nabla_T|\nabla w||^2 & \hbox{for $x\in\{|\nabla w|>0\cap \Omega \}$,} \\
                       0 & \hbox{for $x\in\{|\nabla w|=0\cap \Omega \}$,}
                                                                       \end{array}
                    \right.   \end{eqnarray} 
 where $ \mathcal{\kappa}_i$ are the principal curvatures of the level set of $w$ at $\mathbf x''$ and $\nabla_T$ denotes the orthogonal projection of the gradient along this level set . Setting $w(\mathbf x'')=u(\mathbf x', \mathbf x'')$ and applying this formula together with (\ref{stend}), we finally get (\ref{poincarein}).
 
\hfill $ \Box$

\section{$m$-Liouville theorems for the nonlinear equation}\label{secm}

We now apply  Proposition \ref{liouville}, the linear 0-Liouville theorem, to prove the following  $(d+1)$-Liouville theorem under a strong assumption on the gradient of solutions.

\begin{prop}\label{liou1}
Let $u$ be a monotone solution of (\ref{main}). If there exists $C(n,d)>0$ such that 
\begin{equation}\label{boundthm}
\int_{B_{2R}\setminus B_{R}}  \gamma(\mathbf x')|\nabla_{\mathbf x''} u|^2 d\mathbf x \le C  R^2 g(R),
\end{equation}
for any $g\in \mathcal{G}$.  Then, (\ref{main}) satisfies at most $(d+1)$-Liouville theorem.
\end{prop}

\noindent\textbf{Proof:}   Define $\phi_i(\mathbf x):=\frac{\partial u}{\partial x_i} (\mathbf x)$ for all $i=d+1,\cdots,n$ and $\mathbf x\in\mathbf{R}^n $.  Taking derivative of  (\ref{main}), we get that $\phi_i $ satisfies the following linearized equation $$ -\div (\gamma (\mathbf x') \nabla \phi_i)= \lambda (\mathbf x') f'(u)\phi_i \ \ \text{for all} \ \ \mathbf x \in \mathbf{R}^n.$$
The straightforward calculations show that  
$$ \div (\gamma (\mathbf x') \phi_n^2 \nabla \sigma_i)= 0 \ \ \text{for all} \ \  i=d+1,\cdots,n$$ where $\sigma_i:=\frac{\phi_i}{\phi_n}$.  Note that $\phi_n^2 \sigma_i^2=|\partial_i u |^2$ and from (\ref{boundthm})  for all $i=d+1,\cdots,n$  we have 
  $$
\int_{B_{2R}\setminus B_{R}} \gamma(\mathbf x') \phi_n^2\sigma_i^2 d\mathbf x =\int_{B_{2R}\setminus B_{R}} \gamma(\mathbf x') |\partial_i u |^2  d\mathbf x \le \int_{B_{2R}\setminus B_{R}} \gamma(\mathbf x') |\nabla_{\mathbf x''} u |^2 d\mathbf x \le C R^2 g(R).
$$
Applying Proposition \ref{liouville} with $h(\mathbf x)= \gamma(\mathbf x') \phi_n^2(\mathbf x)$, we get that $(\sigma_i)_{i=d+1}^{n}$ are all constant. Therefore, there exits $(k_i)_{i=d+1}^{n}$ such  that $\sigma_i(\mathbf x)=k_i$ for any  $\mathbf x\in \mathbf R^n$. Clearly $k_n=1$.

From the definition of $\sigma_i$ we get $\frac{\partial u}{\partial x_i}(\mathbf x)=k_i \frac{\partial u}{\partial x_n}(\mathbf x)$ for all $i=d+1,\cdots,n-1$. Therefore, $\nabla _{\mathbf x''} u(\mathbf x)=\frac{\partial u}{\partial x_n} (\mathbf x) (k_{d+1},k_{d+2}, \cdots, k_{n-1}, 1)$. Since $u$ is monotone in $x_n$ direction that is  $\frac{\partial u}{\partial x_n}>0$ we conclude that $\nabla _{\mathbf x''} u(\mathbf x)$ does not change sign for all $x\in\mathbf{R}^n$. Also, note that $u$ is constant along the following directions: $$(\underbrace{0,0,\cdots,0}_{d \ \text{times}},   \underbrace{1,0,\cdots,0,-k_{d+1}}_{s\ \text{times}}), (\underbrace{0,0,\cdots,0}_{d \ \text{times}},\underbrace{0,1,0,\cdots,0,-k_{d+2}}_{s\ \text{times}}),\cdots,  (\underbrace{0,0,\cdots,0}_{d \ \text{times}},\underbrace{0,0,\cdots,0, 1,-k_{n-1}}_{s\ \text{times}})  .$$ 
Therefore, $u$ is a function of $(\mathbf x',\mathbf k\cdot\mathbf  x'')$ where $\mathbf k=(k_{d+1},\cdots,k_{n-1},1)$. 

\hfill $ \Box$

\begin{rem}
Applying the geometric Poincar\'{e} inequality that is given as Theorem \ref{poincare} when $\phi$ is the following standard test function
 $$\phi (x):=\left\{
                      \begin{array}{ll}
                        \frac{1}{2}, & \hbox{if $|x|\le\sqrt{R}$,} \\
                      \frac{ \log R -\log |x|}{{\log R}}, & \hbox{if $\sqrt{R}< |x|< R$,} \\
                       0, & \hbox{if $|x|\ge R$.}
                                                                       \end{array}
                    \right.$$
one can prove Proposition \ref{liou1} for stable solutions as well.  This test function is also used in \cite{bcn,gg2,fg} in order to prove certain results related to the De Giorgi's conjecture.  
\end{rem}

Now we are ready to provide the proof of Theorem  \ref{liou3}.  The idea is to apply the linear 0-Liouville theorem to prove a 0-Liouville theorem for the equation (\ref{main}). 
\\
\\
\noindent\textbf{Proof of Theorem \ref{liou3}:}  Since $u$ is a pointwise stable solution, there exists $v>0$ such that  $$ -\div (\gamma (\mathbf x') \nabla v)= \lambda (\mathbf x') f'(u)v\ \ \text{for all} \ \ \mathbf x \in \mathbf{R}^n.$$
It is straightforward to see that  
\begin{equation}\label{linearv}
 \div (\gamma (\mathbf x') v^2 \nabla \sigma_i)= 0 \ \ \text{for all} \ \  i=d+1,\cdots,n
 \end{equation}
where $\sigma_i:= \frac{\partial u}{\partial_{x_i}}/ {v}$. Therefore, for all  $i=d+1,\cdots,n$ we have $(\sigma_i v)^2\le |\nabla u|^2$ that gives
\begin{equation}\label{liou3sigma}
\int_{B_{R}} \gamma(\mathbf x') v^2 \sigma_i^2\le  \int_{B_R}  \gamma(\mathbf x') |\nabla u|^2.
  \end{equation} 
 To apply Proposition \ref{liouville} we need to find an upper bound for the right-hand side of the  inequality  (\ref{liou3sigma}).   First, we assume that $f$ is a nonnegative nonlinearity.  Multiply both sides of (\ref{main}) with $(u-||u||_{\infty})\phi^2$ where $0\le\phi\le 1 $ is a test function. Since $\lambda (\mathbf x')f(u)(u-||u||_{\infty})\le 0$,  we have 
\begin{equation}\label{supersol}
\hfill -\div(\gamma(\mathbf x') \nabla u) (u-||u ||_\infty)\phi^2 \le 0  \ \ \text{in}\ \ \ \mathbf{R}^n=\mathbf{R}^d\times\mathbf{R}^s.
  \end{equation}
On the other hand, for the case $tf(t)\le 0$ a similar differential inequality holds. Note that  multiplying  both sides of (\ref{main}) with $u\phi^2$ we have   
\begin{equation}\label{supersoluf}
\hfill -\div(\gamma(\mathbf x') \nabla u) u \phi^2\le 0  \ \ \text{in}\ \ \ \mathbf{R}^n=\mathbf{R}^d\times\mathbf{R}^s.
  \end{equation}
Now, integrating both sides of (\ref{supersol}) and (\ref{supersoluf}) and using the fact that $u$ is bounded we obtain
\begin{eqnarray*}
\int_{\mathbf{R}^n} \gamma(\mathbf x')  |\nabla  u|^2\phi^2&\le& k \int_{\mathbf{R}^n} \gamma(\mathbf x')  |\nabla u ||\nabla \phi |\phi \\&\le&  k \left(  \int_{\mathbf{R}^n} \gamma (\mathbf x') |\nabla u|^2 \phi^2 \right)^{\frac{1}{2}} \left(  \int_{\mathbf{R}^n} \gamma(\mathbf x') |\nabla \phi|^2  \right)^{\frac{1}{2}} 
  \end{eqnarray*}
  where $k$ is a constant that  only  depends  on $||u||_{\infty}$.  Set the test function $\phi$ to be the standard test function that is $\phi=1$ in $B_R$ and $\phi=0$ in $\mathbf{R}^n\setminus B_{2R}$ where $||\nabla \phi||_{L^\infty(B_{2R}\setminus B_R)}\le k R^{-1}$. Then, we have  
  \begin{eqnarray}\label{decay}
\nonumber\int_{B_{R}} \gamma(\mathbf x')  |\nabla  u|^2 &\le& k R^{-2} \int_{B_{2R}\setminus B_{R}}\gamma (\mathbf x') d\mathbf x\\&\le& k R^{s-2} \int_{B_{2R}} \gamma (\mathbf x') d\mathbf x'.
    \end{eqnarray}
 Let (\ref{liou3g}) hold in dimensions $n\le d+4$, then  $n-d-2=s-2\le 2$ and  
 $$R^{s-2} \int_{B_{2R}} \gamma (\mathbf x') d\mathbf x' \le  R^2 \int_{B_{2R}} \gamma (\mathbf x') d\mathbf x' \le k R^2 g(R).$$
Similarly, if   (\ref{liou3rg}) holds when $n\le d+3$ then $n-d-2=s-2\le 1$ and 
  $$R^{s-2} \int_{B_{2R}} \gamma (\mathbf x') d\mathbf x'\le  R \int_{B_{2R}} \gamma (\mathbf x') d\mathbf x' \le k R^2 g(R).$$
 Therefore, from (\ref{decay}) and (\ref{liou3sigma}) we get
 $$ \int_{B_{R}} \gamma(\mathbf x') v^2 \sigma_i^2 \le k R^2 g(R) \ \ \text{for all} \ \ i=d+1,\cdots,n.$$
 Set $h(\mathbf x)=\gamma(\mathbf x')v^2$ and $\sigma=\sigma_i$ in  Proposition \ref{liouville} for all $ i=d+1,\cdots,n$ to obtain that all $\sigma_i$ are constant.  By similar discussions as in the proof of Proposition \ref{liou1} we have $u(\mathbf x',\mathbf x'')=w(\mathbf x',\mathbf k\cdot \mathbf x'')$ such that $\mathbf k\in\mathbf R^s$ and  $| \mathbf k|=1$. Note that $w$ satisfies 
\begin{equation}
\label{supersolw}
\hfill  (w-||w ||_\infty)  \div\left(\gamma(\mathbf x') \nabla (w-||w ||_\infty) \right) \ge 0  \ \ \text{in}\ \ \ \mathbf{R}^{d+1}=\mathbf{R}^d\times\mathbf{R},
  \end{equation}
 where $f(w)\ge 0$ and similarly
 \begin{equation}
\label{supersolwuf}
\hfill  w \div\left(\gamma(\mathbf x') \nabla w \right) \ge 0  \ \ \text{in}\ \ \ \mathbf{R}^{d+1}=\mathbf{R}^d\times\mathbf{R},
  \end{equation}
 where $wf(w)\le 0$. The fact that  $w$ is bounded and  satisfies either (\ref{supersolw}) or (\ref{supersolwuf}) in dimension $d+1$ and decay estimates (\ref{liou3g}) and (\ref{liou3rg}) hold for $\gamma$ imply that  
 \begin{eqnarray*}
\int_{B_{R}} \gamma(\mathbf x')    (w-||w||_\infty)^2 d\mathbf x \le k R \int_{B_{R}} \gamma(\mathbf x') d\mathbf x'\le  R^2 g(R)\\
\int_{B_{R}} \gamma(\mathbf x')    w ^2 d\mathbf x \le k R \int_{B_{R}} \gamma(\mathbf x') d\mathbf x'\le  R^2 g(R)
  \end{eqnarray*}
  where $k$ is a positive constant independent of $R$. Hence applying Proposition \ref{liouville} again for (\ref{supersolw}) and (\ref{supersolwuf}) we obtain that $w$ is constant. Therefore, $u$ is constant.

\hfill $ \Box$

Note that to apply Proposition \ref{liou1} one needs to  have a $L^2(B_R)$ upper bound on $|\nabla u|$ that we call the energy bound. In  what follows we give such an energy bound in terms of weight functions $\lambda$ and $\gamma$. The following lemma holds for subsolutions of (\ref{main}) as well. By subsolution we  mean the inequality $``\le "$ holds in (\ref{main}).

\begin{lemma}\label{enerbound}
Let  $u$ be a bounded solution of (\ref{main}) with any $f\in C^1(\mathbf R)$. Then 
\begin{equation}\label{boundf}
\int_{ B_{R}} \gamma(\mathbf x')  |\nabla u|^2 \ d\mathbf x \le k R^{s} \int_{ B_{R}} \left\{ \lambda(\mathbf x')+R^{-2} \gamma( \mathbf x') \right\} d\mathbf x',
\end{equation}
where the positive constant $k$ is independent of $R$.
\end{lemma}
\noindent\textbf{Proof:} Multiply both sides of (\ref{main}) with $(||u||_\infty+u)\phi^2$ when $0\le\phi\le 1$ is a test fucntion. Then, integrating by parts we get
 \begin{equation*}
 \int_{\mathbf R^n} \gamma(\mathbf x') \nabla u\cdot \nabla \left( \phi^2 (||u||_\infty+u) \right) \le \int_{\mathbf R^n}  \lambda(\mathbf x') f(u)  (||u||_\infty+u)  \phi^2.
 \end{equation*}
 Simplifying this inequality and keeping the square of gradient of $u$ in the left hand side, we end up with
  \begin{equation}\label{grad2}
 \int_{\mathbf R^n} \gamma(\mathbf x') |\nabla u|^2 \phi^2 \le \int_{\mathbf R^n}  \lambda(\mathbf x') f(u)  (||u||_\infty+u)  \phi^2+4||u||_{\infty} \int_{\mathbf R^n} \gamma(\mathbf x')  |\nabla u| |\nabla\phi| \phi  
 \end{equation}
 We now define the positive  constants $k$ and $\epsilon$ such that $ 2||f(u)||_\infty ||u||_\infty<k<\infty$ and $0<\epsilon< (4||u||_\infty)^{-1}$. Applying the Young's inequality\footnote{For any positive $\epsilon$ and any $a,b\in\mathbf{R}$, $ab\le \epsilon a^2+\frac{1}{4\epsilon} b^2$.} for the last term in right hand side of (\ref{grad2}) we get
 \begin{equation}\label{lems}
 \left(1- 4||u||_\infty\epsilon \right)  \int_{\mathbf R^n} \gamma(\mathbf x') |\nabla u|^2 \phi^2 \le  k \int_{\mathbf R^n}  \lambda(\mathbf x') \phi^2+ ||u||_\infty \epsilon^{-1}  \int_{\mathbf R^n}  \gamma(\mathbf x') |\nabla \phi|^2.
 \end{equation}
Finally, set $\phi$ to be the standard smooth test function that is $\phi=1$ in $B_R$ and $\phi=0$ in $\mathbf{R}^n\setminus B_{2R}$ with $||\nabla\phi||_{L^\infty({B_{2R}})}< k R^{-1}$. This proves (\ref{boundf}).

\hfill $ \Box$

In the statement of Lemma \ref{enerbound}, there is no assumption on the monotonicity of the solutions. However, monotonicity is a crucial assumption to derive $m$-Liouvile theorems when $m\ge 1$.   In other words, assuming the monotonicity of solutions we get a stronger upper bound on the energy of solutions. Before we discuss the new upper bound on the energy $E_R$, let us mention that applying some standard elliptic estimates to bounded solutions of (\ref{main}) gives us $|\nabla u|\in L^{\infty}(\mathbf{R}^n)$.  Indeed, assume that $u$ is a bounded solution of either (\ref{main}) when $\frac{|\nabla \gamma|}{\gamma}, \frac{\lambda}{\gamma}\in L^{\infty}(\mathbf{R}^n)$ or equivalently (\ref{vector}) when $\mathbf a,b\in L^{\infty}(\mathbf{R}^n)$.  Then applying interior $W^{2,p}$ estimates with $p>n$ to $-\Delta u+\mathbf{a}(\mathbf x') \cdot \nabla u=b(\mathbf x') f(u)\in L^{\infty}(B_2(\mathbf y))$ for every $\mathbf y\in\mathbf{R}^n$, we get 
 \begin{equation*}
||u||_{W^{2,p}(B_1(\mathbf y))} \le k \{||u||_{L^\infty(B_2(\mathbf y)) }+ ||f(u)||_{L^p((B_2(\mathbf y))}  \} \le k,
\end{equation*}
where $k$ is independent of $\mathbf y$.    Using the Sobolev embedding $W^{2,p}(B_1(\mathbf y)) \subset C^1\left(\overline{B_1(\mathbf y)}\right)$  for $p>n$ and any $ \mathbf y\in\mathbf{R}^n$, we have $u\in C^1(\mathbf{R}^n)$ and $|\nabla u|\in L^{\infty}(\mathbf{R}^n). $

\hfill $ \Box$

\begin{lemma}\label{boundinf}
Let $u$ be a bounded monotone solution of (\ref{main}) for any $f\in C^1(\mathbf{R})$ and $$\lim_{x_n\to \infty} u(\mathbf x''',x_n)=1\ \ \ \text{for all}\ \  \mathbf x=(\mathbf x''', x_n)\in\mathbf{R}^{n-1}\times\mathbf{R}.$$ 
Then 
\begin{equation}\label{energybound}
E_R(u) \le k \int_{\partial B_R} \gamma(\mathbf x') dS(\mathbf x),
\end{equation}
where  the positive constant $k$ is independent of $R$ and $E_R(u)$ is the energy functional defined by
$$E_R(u):=\frac{1}{2}\int_{B_{R}} \gamma(\mathbf x')|\nabla u|^2 d\mathbf x - \int_{B_{R}} \lambda ( \mathbf x') (F(u)-F(1)) d \mathbf x. $$
\end{lemma}

\noindent\textbf{Proof:}  Define $u^t(\mathbf x)=u(\mathbf x''',x_n+t)$ for  $t\in\mathbf{R}$. Note that the shifted function $u^t$ satisfies (\ref{main}) that is 
\begin{eqnarray}
\label{maint}
-\div(\gamma(\mathbf x') \nabla u^t) =\lambda(\mathbf x') f(u^t)   \quad \text{in}\ \ \mathbf{R}^n. 
  \end{eqnarray}
  Moreover, the following monotonicity and decay conditions hold.
  \begin{eqnarray}
\label{monot}
 \left\{\begin{array}{lcl}
  \partial_t u^t(\mathbf x) >0 \ \ \text{in}\ \ \mathbf{R}^n\\ 
 \lim_{t\to\infty} u^t(\mathbf x)=1\ \ \text{in}\ \ \mathbf{R}^n\\
 |\nabla u^t| \in L^{\infty}(\mathbf{R}^n).  
\end{array}\right.
  \end{eqnarray} 
\noindent Step 1: We claim that the following decay estimate holds for any $R>1$
\begin{equation}\label{decayclaim}
\lim_{t\to\infty} E_R(u^t)=0.
\end{equation}
To prove (\ref{decayclaim}) we apply the properties of $u^t$ given in (\ref{monot}). Since   $\lim_{t\to\infty} u^t(\mathbf x)=1$ for any  $\mathbf x\in\mathbf{R}^n$,  for any $R>1$ we get 
$$\lim_{t\to\infty} \int_{B_{R}}\lambda(\mathbf x') (F(u^t) -F(1))  = 0.$$
From definition of the energy functional  $E_R$ we only need to prove that for any $R>1$
\begin{equation}\label{decayclaim2}
\lim_{t\to\infty}  \int_{B_{R}} \gamma( \mathbf x') |\nabla u^t|^2  \to 0 \ \ \text{ for any } R>1.
  \end{equation}
Multiply both sides of (\ref{maint}) with $u^t-1$ and do  integration by parts on $B_R$ to end up with 
\begin{equation*}
\label{}
 \int_{B_{R}} \gamma(\mathbf x') |\nabla u^t|^2=\int_{\partial B_{R}}\gamma(\mathbf x') (u^t-1)\partial_\nu u^t +\int_{B_{R}}  \lambda(\mathbf x') f(u^t) (u^t-1). 
  \end{equation*}
Taking the limit of both sides as $t\to\infty$ and using again the fact that $\lim_{t\to\infty} u^t(\mathbf x)=1$ for any  $\mathbf x\in\mathbf{R}^n$, one can get (\ref{decayclaim2}). This finishes the proof of (\ref{decayclaim}). 

\noindent Step 2:  The following upper bound holds for the energy of $u$
\begin{equation}\label{energyt}
E_R(u) \le  E_R(u^t) +k \int_{\partial B_{R}} \gamma(\mathbf x') dS(\mathbf x) \ \ \text{for all}\ \ t \in\mathbf{R^+},
\end{equation}
where $k$ is a contact and independent of $R$.     Differentiating the energy functional of  $u^t$ gives us 
\begin{eqnarray}\label{der-energy}
\partial_t E_R(u^t)= \int_{B_{R}} \gamma(\mathbf x') \nabla u^t\cdot \nabla (\partial _t u^t)- \int_{B_{R}}\lambda(\mathbf x') f(u^t) \partial_t u^t.
    \end{eqnarray}
Now, multiply (\ref{maint}) with $\partial_t u^t$ and perform integration by parts on $B_R$ to get 
  \begin{eqnarray}
\label{partialt}
\hfill  \int_{B_{R}} \gamma(\mathbf x') \nabla u^t\cdot \nabla (\partial _t u^t)- \int_{\partial B_{R}} \gamma(\mathbf x') \partial_\nu u^t \partial_t u^t &=&\int_{B_{R}} \lambda(\mathbf x') f(u^t) \partial_t u^t.
  \end{eqnarray}
  Note that the integral terms   $\int_{B_{R}} \lambda(\mathbf x') f(u^t) \partial_t u^t$ and $\int_{B_{R}} \gamma(\mathbf x') \nabla u^t\cdot \nabla (\partial _t u^t)$  are common in (\ref{der-energy}) and (\ref{partialt}).  So,  combining these two integral equalities we get  a simplified form for the derivative of the energy of $u^t$ 
  \begin{eqnarray}\label{der-energy-bd}
\partial_t E_R(u^t)= \int_{\partial B_{R}} \gamma(\mathbf x') \partial_{\mathbf\nu} u^t \partial_t u^t .
   \end{eqnarray}
 Note that the directional derivative of $u^t$ is $\partial_{\mathbf\nu} u^t(\mathbf x)=\mathbf\nu(\mathbf x) \cdot\nabla u^t(\mathbf x)=\mathbf\nu(\mathbf x)\cdot\nabla u(\mathbf x''',x_n+t)$ when $||\mathbf\nu||=1$.  Therefore, $-||\nabla u||_{L^\infty(\mathbf R^n)} \le \partial_{\mathbf\nu} u^t (\mathbf x)\le ||\nabla u||_{L^\infty(\mathbf R^n)} $. From this and the fact that $\partial_t u^t(\mathbf x)>0$ for all $\mathbf x\in\mathbf R^n$ and $t\in\mathbf R^+$, we get 
  \begin{equation}\label{der-energy-bdM}
\partial_t E_R(u^t) \ge -  ||\nabla u||_{L^\infty(\mathbf R^n)}  \int_{\partial B_{R}}  \gamma(\mathbf x') \partial_t u^t  d S(\mathbf x).
   \end{equation}
On the other hand, basic integration shows that 
 \begin{equation*}
 E_R(u) = E_R(u^t)- \int_{0}^{t} \partial_s E_R(u^s) ds. 
  \end{equation*}
From this and (\ref{der-energy-bdM}) we get 
 \begin{equation*}
 E_R(u)\le E_R(u^t) +||\nabla u||_{L^\infty(\mathbf R^n)}  \int_{0}^{t}  \int_{\partial B_{R}} \gamma(\mathbf x') \partial_s u^s dS(\mathbf x) ds. 
   \end{equation*}
Therefore,
 $$ E_R(u)\le E_R(u^t) +||\nabla u||_{L^\infty(\mathbf R^n)}  \int_{\partial B_{R}} \gamma(\mathbf x') (u^t-u) dS(\mathbf x).$$
Note that from the definition of $u^t$, we get $u(\mathbf x)<u^t(\mathbf x)$ for all $\mathbf x\in\mathbf{R}^n$ and $t\in \mathbf{R^+}$ and then $0<u^t(\mathbf x)- u(\mathbf x)<||u||_{L^\infty(\mathbf R^n)}$. Set $k=||\nabla u||_{L^\infty(\mathbf R^n)} ||u||_{L^\infty(\mathbf R^n)}$, this finishes the proof of (\ref{energyt}). To complete the proof, just take the limit of (\ref{energyt}) as $t\to\infty$ in the light of (\ref{decayclaim}). 
 
\hfill $ \Box$

Now, we prove an elementary inequality that compares the surface integral with the volume integral.   
\begin{lemma}\label{surface}
Let $s \geq 2$, $d\ge 1$ and $\gamma \in C^\infty (\mathbf R^d)$ be positive. Then
\[
\int_{\partial B_R}\gamma(\mathbf x') dS (\mathbf x) \leq k R^{s-1} \int_{B_R}\gamma(\mathbf x') d\mathbf  x'
\]
where $k$ is independent of $R$.
\end{lemma}
\noindent\textbf{Proof:}   For a general surface $x_n = \phi(\mathbf x''')$,
the surface area element is $dA = \sqrt{1+ |D\phi|^2} dx_1 \cdots dx_{n-1}$.
For the sphere $\phi(\mathbf x''')=(R^2 - |x_1|^2-|x_2|^2-\cdots |x_{n-1}^2|)^{1/2}$ and therefore
\[
dA=\sqrt{1+ |D\phi|^2} dx_1 \cdots dx_{n-1} = \frac{R}{\phi}dx_1\cdots dx_{n-1}.
\]
Integrating out the $\mathbf x''$-variable, we have
\begin{equation*}
\int_{\partial B_R}\gamma(\mathbf x') dS (\mathbf x) = \int_{B_R}\gamma(\mathbf x') w(R,\mathbf x') d\mathbf x'
\end{equation*}
for some weight function $w(R,\mathbf x')\geq 0$. We now prove that
\begin{equation}\label{surfacequ}
w(R,\mathbf x')=k_s R(R^2-|\mathbf x'|^2)^{\frac{s-2}{2}},
\end{equation}
where $k_s$ is a constant independant of $R$. This proves the lemma since $w(R,\mathbf x')\leq k_s R^{s-1}$ whenever $s\geq 2$.    Rewrite $\phi=(\rho^2-|\mathbf y|^2)^{1/2}$, where $\rho^2=R^2-|\mathbf x'|^2$ and $\mathbf x''=(\mathbf y, x_n)\in R^s$. The weight  function is then
\begin{equation}\label{surfacequ2}
w(R, \mathbf x')=\int_{|\mathbf y|<\rho}\frac{R}{\phi} d\mathbf y=R\int_{|\mathbf y|<\rho}\frac{d\mathbf y}{(\rho^2-|\mathbf y|^2)^{1/2}} = k_s R \rho^{s-2} 
\end{equation}
where $k_s:=\int_{B_1^{s-1}}  \frac{d\mathbf z}{   (1-|\mathbf z|^2)^{1/2}  } $ and $  B_1^{s-1}$ is the unit ball in $\mathbf R^{s-1}$. From the definition of $\rho$, this proves  (\ref{surfacequ}).

\hfill $ \Box$

We are now ready to see the proof of Theorem \ref{liou2}. 
\\
\\
\noindent\textbf{Proof of Theorem \ref{liou2}:} Without loss of generality we assume that $F(-1)\ge F(1)$. Therefore, from the assumptions $F(u)-F(1)\le 0$ that gives us
$$ \int_{B_{R}} \lambda ( \mathbf x') (F(u)-F(1)) d \mathbf x\le0.$$
From this and Lemma \ref{boundinf} we get the following bound on the gradient of solutions
\begin{equation*}
\int_{B_{R}} \gamma(\mathbf x')|\nabla u|^2 d\mathbf x \le k \int_{\partial B_R} \gamma(\mathbf x') dS(\mathbf x),
\end{equation*}
where $k$ is a constant independent of $R$. Applying Lemma \ref{surface} we change the upper bound to a volume integral of $\gamma$ that  is 
\begin{equation}\label{thm2es}
\int_{B_{R}} \gamma(\mathbf x')|\nabla u|^2 d\mathbf x \le k  R^{s-1} \int_{B_R}\gamma(\mathbf x') d\mathbf x'.
\end{equation}
The key point is to apply Proposition \ref{liou1} to show that solutions are at most $(d+1)$-dimensional. So, we need to make sure that (\ref{boundthm}) holds.   From (\ref{thm2es}) we only need 
\begin{equation}\label{thm2es2}
  \int_{B_R}\gamma(\mathbf x') d \mathbf x'  \le k  R^{3-s} g(R),
\end{equation}
for any $g\in \mathcal{G}$.    Note that for a positive $\gamma$ to satisfy (\ref{thm2es2}) we need to assume that $s\le 3$ and also we assumed $s\ge 2$ to prove Lemma \ref{surface}. Therefore, for $s=2$ that is $n=d+2$ we assume that 
\begin{equation*}
  \int_{B_R}\gamma(\mathbf x') d \mathbf x'  \le k  R g(R),
\end{equation*}
and for $s=3$ that is $n=d+3$ we assume that 
\begin{equation*}
  \int_{B_R}\gamma(\mathbf x') d \mathbf x'  \le k  g(R). 
\end{equation*}
This finishes the proof for the case $F(-1) \ge F(1)$. Note that if $F(-1) < F(1)$, replace $u(\mathbf x''',x_n)$ with $-u(\mathbf x''',-x_n)$ and apply the same argument. 

\hfill $ \Box$

\emph{Acknowledgement}. I am enormously grateful to the anonymous referee for carefully reading the paper and providing constructive comments.

      \end{document}